\documentclass[12pt]{amsart}
\usepackage{amsmath}

\newtheorem{thm}{Theorem}
\newtheorem{lem}[thm]{Lemma}
\newtheorem{prop}[thm]{Proposition}

\theoremstyle{remark}
\newtheorem{remark}{Remark}[section] 

\newcommand{\Z}{{\mathbf Z}}
\newcommand{\Q}{{\mathbf Q}}

\newcommand{\R}{{\mathbf R}}
\newcommand{\C}{{\mathbf C}}

\newcommand{\trace}{\operatorname{tr}}

\newcommand{\OO}{\mathfrak{O}}

\newcommand{\weil}{U}
\newcommand{\weilp}{U_p}
\newcommand{\torus}{{\mathcal C}_A}

\newcommand{\PP}{P_{\nu}}
\newcommand{\norm}[1]{\left\| #1 \right\|_2}
\newcommand{\supnorm}[1]{\left\| #1 \right\|_\infty}

\newcommand{\Exp}{E_{\nu}}
\newcommand{\cond}{\operatorname{cond}}

\newcommand{\CA}{F}

\newcommand{\HN}{H_N}
\newcommand{\HP}{H_p}
\newcommand{\UN}{U_N}

\newcommand{\fp}{\mathbb{F}_p}
\newcommand{\fpt}{\mathbb{F}_{p^2}}
\newcommand{\fpn}{\mathbb{F}_{p^2}^1}

\title[Supremum norms for Hecke eigenfunctions]{Bounds on supremum norms for Hecke eigenfunctions of quantized cat maps}

\author{P\"ar Kurlberg}

\email{kurlberg@math.kth.se}
\address{Department of Mathematics\\ 
Royal Institute of Technology\\
SE-100 44 Stockholm  \\
Sweden}

\date{April 28, 2006}

\thanks{Partially supported by a grant from 
the G\"oran Gustafsson Foundation, 
the National Science
  Foundation (DMS-0071503), the Royal Swedish Academy of Sciences, and
  the Swedish Research Council.  } 

\begin{document}
\begin{abstract}
  We study extreme values of desymmetrized eigenfunctions (so called
  Hecke eigenfunctions) for the quantized cat map, a quantization
  of a hyperbolic linear map of the torus.
  In a previous paper it was shown that for prime values of the
  inverse Planck's constant $N=1/h$, such that the map is
  diagonalizable (but not upper triangular) modulo
  $N$, 
the Hecke eigenfunctions are uniformly bounded.  The purpose
  of this paper is to show that the same holds for any prime $N$
  provided that the map is not upper triangular modulo $N$.
   We also find that the supremum norms
  of Hecke eigenfunctions are $\ll_\epsilon N^\epsilon$ for all
  $\epsilon>0$ in the case of $N$ square
  free.

\end{abstract}

\maketitle

\section{Introduction}

The behavior of eigenfunctions, such as their value distribution and
extreme values, of classically chaotic quantum systems has received
considerable attention in the past few years
\cite{berry-regular-irregular-semiclassical,
  hejhal-rackner-gaussian-value-distribution, sarnak-arithmetic-qc,
  aurich-steiner-statistical-properties-highly-excited,
  iwaniec-sarnak-supnorms,
  aurich-backer-schubert-taglieber-maximum-norms} (see
section~\ref{subsec:comparison} for a discussion of these and related
results.)
The aim of this note is to investigate supremum norms for
eigenfunctions in the context of ``quantized cat maps''.  The
classical maps, known as ``cat maps'', are given by the action of
hyperbolic linear maps $A \in SL_2(\Z)$ on the two dimensional torus.
A procedure to 
quantize such 
maps was first introduced by Berry and Hannay in \cite{BH}, and was further
developed in \cite{Knabe,DE,DEGI,klimek-leniewski-maitra-rubin,Zelditch97,
  cat1,gurevich-hadani-hanny-berry-equivariant,
  mezzadri-multiplicative-cat-maps}.
The quantization procedure restricts Planck's constant to be of the
form $h = 1/N$ where $N$ is a positive integer.  For each integer $N\geq
1$, let $\UN(A)$ denote the quantization of $A$ as a unitary operator
on the Hilbert space of states $\HN$.  (For more
details, see section~\ref{sec:background} and
references mentioned therein.)   For certain values of $N$,
$\UN(A)$ can have very large spectral degeneracies, and in
\cite{cat1} it was observed that these degeneracies are connected to
quantum symmetries of $\UN(A)$, namely a commutative group of unitary
operators which contains $\UN(A)$.  These operators are called Hecke
operators in analogy with the classical theory of modular forms.  The
eigenspaces of $\UN(A)$ admit an orthonormal basis consisting of
eigenfunctions of all the Hecke operators, which are called {\em Hecke
  eigenfunctions}.
%
In \cite{catsup} it was shown that for prime values of $N$ for which
the map is diagonalizable (but not upper triangular) modulo $N$ (the
so called ``split primes''), normalized Hecke eigenfunctions are
uniformly bounded, and, by using results of N. Katz \cite{katz-kr-sums},
that their value distribution is given by the semi-circle law.
The aim of this paper is to show that the same bound 
holds for all odd prime $N$, provided that the
map is not upper triangular modulo $N$:
%
%
\begin{thm}  
\label{thm:main}
If $N$ is an odd inert prime, $\psi \in \HN$ is a Hecke eigenfunction
normalized so that $\norm{\psi}=1$,
and $A$ is not upper triangular modulo $N$,
then
$$
\supnorm{\psi} \leq \frac{2}{\sqrt{1+1/N}}
$$
\end{thm}
{\em Remark: }  
By a different method (namely expressing the Hecke eigenfunctions in
terms of perverse sheaves and then bounding the dimensions of certain
cohomology groups) the same result has independently been obtained by
Gurevich and Hadani in \cite{gurevich-hadani-hecke-sup-norms}.

Theorem~\ref{thm:main} together with Theorem~2 of \cite{catsup}, and
a short treatment of the case of ramified primes gives that Hecke
eigenfunctions are 
uniformly bounded for all prime values of $N$ as long as $A$ is not
upper triangular modulo $N$.
\begin{thm}  
\label{thm:summary}
If $N$ is an odd prime, $\psi \in \HN$ is an $L^2$-normalized Hecke
eigenfunction, and $A$ is not upper triangular modulo $N$, then 
$$
\supnorm{\psi} \leq 2
$$
On the other hand, if $A$ is upper triangular modulo $N$, there exists
Hecke eigenfunctions $\psi$ such that
$$
\supnorm{\psi} \geq \sqrt{N/2}
$$
\end{thm}
{\em Remark: } Since $A$ is assumed to be hyperbolic, it can only be
upper triangular modulo $N$ for finitely many $N$.

For composite $N$, it was shown in \cite{catsup} that $\supnorm{\phi}
\ll_\epsilon N^{3/8+\epsilon}$ for all
$\epsilon>0$ if $\phi \in \HN$ is a Hecke
eigenfunction.  However, for $N$ square free, this bound can be
improved considerably by using Theorem~\ref{thm:summary}.
\begin{thm}
\label{thm:square-free-bound}
If $N$ is square free and $\phi \in \HN$ is an $L^2$-normalized Hecke
eigenfunction, then 
$$
\supnorm{\phi} \ll_\epsilon N^\epsilon
$$  
for all $\epsilon>0$.
\end{thm}
However, it should be noted that this bound does not hold for general
integers $N$ --- 
Olofsson has shown \cite{olofsson-big-catmap-supnorm} that there
exists a subsequence $\{N_i\}_{i=1}^\infty$ for which there exists
Hecke eigenfunctions $\phi_i \in H_{N_i}$ satisfying $\supnorm{\phi_i}
\geq N_i^{1/4}$.

\subsection{Comparison with eigenfunctions of the Laplacian}
\label{subsec:comparison}




A rich class of examples of chaotic classical dynamics are
given by the geodesic flow on negatively curved compact Riemannian
manifolds.  In this setting, the quantization of the dynamics is
essentially given by the Laplace-Beltrami operator $\Delta$ acting on
smooth functions on the manifold, and the eigenfunctions of the
quantized cat map can be viewed as an analogue of eigenfunctions of
$-\hbar^2 \Delta$, i.e., functions $\psi_\lambda$ (always assumed to
be normalized to have $L^2$-norm one) satisfying
$
-\hbar^2 \Delta \psi_\lambda   = \lambda \psi_\lambda
$
where $\hbar$ is Planck's constant and $\lambda$ is the energy
associated with the eigenstate $\psi_\lambda$.  Keeping the energy
fixed and letting $\hbar \to 0$ (the semiclassical limit) is
equivalent to setting $\hbar=1$ and letting $\lambda \to \infty$ (the
large energy limit), and when making comparisons with the cat map we
should think of $N$, the inverse Planck's constant, to be of size
$\sqrt{\lambda}$.  

As a model for eigenfunctions in the case of classical chaotic
dynamics, Berry proposed \cite{berry-regular-irregular-semiclassical}
a superposition of random plane waves.  Consequently, eigenfunctions
should have a Gaussian value distribution, and this prediction is
matched very well by numerics for certain arithmetic surfaces
\cite{hejhal-rackner-gaussian-value-distribution,aurich-steiner-statistical-properties-highly-excited}.
As for extreme values, the random wave model predicts (see
\cite{hejhal-rackner-gaussian-value-distribution,
  sarnak-arithmetic-qc}, and also
\cite{aurich-backer-schubert-taglieber-maximum-norms} for numerics
specifically investigating large values of eigenfunctions) 
that $\supnorm{\psi_\lambda}$ should be of order $\sqrt{\log
  \lambda}$.
As for rigorous results, a well known bound 
(e.g., see \cite{sogge-sup-norm})
valid for any compact
Riemannian manifold  of dimension $n$ is that
$
\supnorm{\psi_\lambda} \ll \lambda^\frac{n-1}{4},
$
and
the bound is sharp as
can be seen by considering zonal harmonics on the sphere.
(For cat maps, the corresponding bound $\supnorm{\psi} \ll \sqrt{N}$ is
trivial.)
It is of considerable interest to improve this bound using various
dynamical properties of the geodesic flow.  For real analytic
{\em surfaces}, Sogge and Zelditch proved 
\cite{sogge-zelditch-max-eigenfunction-growth} that if the flow is
ergodic, then $\supnorm{\psi_\lambda} = o(\lambda^{\frac{n-1}{4}})$.
Moreover, in the case of negative
curvature, the slightly stronger bound $ \supnorm{\psi_\lambda} \ll
\lambda^\frac{n-1}{4}/\log \lambda $ follows from 
Berard \cite{berard}, but this is probably quite far from the truth,
especially in dimension two.
In fact, Iwaniec and Sarnak conjectured \cite{iwaniec-sarnak-supnorms}
that for surfaces of constant negative curvature,
$
\supnorm{\psi_\lambda} \ll_\epsilon \lambda^\epsilon
$
holds for all $\epsilon>0$.
In this direction, they also 
proved that for certain arithmetic quotients of the upper half plane,
$
\supnorm{\psi_\lambda} \ll_\epsilon \lambda^{5/24 +\epsilon},
$
for all $\epsilon>0$, as well as  that the lower bound
$
\supnorm{\psi_\lambda} \gg \sqrt{\log \log \lambda}
$
holds for infinitely many eigenvalues.  (To be precise, they assume
that $\{\psi_\lambda\}$ is a basis of Hecke eigenfunctions.)
Note that  their results, as well as their conjecture, 
are consistent
with the random wave model prediction for extreme values of
eigenfunctions.  However, it should also be noted that in higher
dimensions, 
eigenfunctions can have rather large supremum norms even
though the curvature is negative --- in
\cite{rudnick-sarnak-eigenstates-arithmetic-hyperbolic-mflds}, 
Rudnick and Sarnak showed that for
certain arithmetic manifolds of dimension three,
$\supnorm{\psi_{\lambda}} \gg \lambda^\frac{1}{4}$ for infinitely many
eigenvalues.  

In the integrable case eigenfunctions are better understood, and it is
sometimes possible to infer fairly refined information about the
geometry of the manifold from the growth of eigenfunctions.
For irrational flat tori eigenfunctions are uniformly bounded, and in
\cite{toth-zelditch-bounded-eigenfunctions,toth-zelditch-revisited}
Toth and Zelditch proved a partial converse: under certain assumptions
(bounded eigenvalue multiplicity, ``complete quantum integrability''
and ``bounded complexity'') it turns out that
uniformly bounded eigenfunctions implies that the metric is flat.
(However, note that rational flat tori have unbounded
multiplicities and hence unbounded eigenfunctions.)
Moreover, Sogge and Zelditch has shown
\cite{sogge-zelditch-max-eigenfunction-growth} (also see
\cite{toth-zelditch-revisited}) that for manifolds with completely
integrable flow satisfying a certain non-degeneracy condititon,
non-flatness is equivalent with large growth rates of the $L^\infty$
and $L^p$-norms for a subsequence of eigenfunctions.  Further, they
also showed that if $M$ is manifold with an infinite subsequence of
maximal growth eigenfunctions (i.e., $\supnorm{\phi_\lambda} \gg
\lambda^{(n-1)/4}$) then there exists
a point $x \in M$ for which the set of directions of geodesic loops at
$x$ has positive measure.  In particular, if $M$ is a real analytic
surface with maximal eigenfunction growth, then $M$ is topologically
either a $2$-sphere or the real projective plane.

For further reading, see
\cite{sarnak-arithmetic-qc,
jakobson-survey,sarnak-spectra-of-hyperbolic-surfaces,zelditch-quantum-ergodicity-and-mixing}
for some nice surveys.

{\bf Acknowledgements:} I would like to thank Zeev Rudnick for 
helpful discussions.

\section{Background}
\label{sec:background}

\subsection{Notation}
Since we only will be concerned with Planck's constant $h$ taking
values among inverse primes, we will use the notation $h=1/p$ where
$p$ is a prime (rather than $h=1/N$.)
$\fp$ will denote the finite field with $p$
elements, $\fpt$ the the finite field with $p^2$ elements, and 
$$
\fpn :=
\{ x \in \fpt : N_{\fp}^{\fpt}(x) = 1\}
$$ 
is the kernel of the norm map.  Further, $\psi : \fp \to \C^\times$
will denote a nontrivial additive character of $\fp$.

\subsection{Classical dynamics}
The classical dynamics are given by a hyperbolic linear map $A \in
SL_2(\Z)$ acting on the phase space $T^2=\R^2/\Z^2$, and the time
evolution on a classical observable $f \in C^\infty(T^2)$ is given by
$f \to f \circ A$.  
Since $A$ is assumed to be hyperbolic, the eigenvalues of $A$ are 
powers of a fundemental unit in a real quadratic field $K$.  
For simplicity, we will assume that $A \equiv I
\mod 2$, but we will outline how this restriction can be avoided in
Remark~\ref{rem:no-parity}.

\subsection{The Hilbert space of states}
The Hilbert space of states
may be identified with
$\HP =L^2(\fp)$, where the inner product is given by
\begin{equation}
  \label{eq:inner-product}
\langle f,g \rangle :=
\frac{1}{p}\sum_{x=1}^p f(x) \overline{g(x)}.
\end{equation}

\subsection{The quantum propagator and the Weil representation}
\label{sec:weil-representation}
We will use the quantization procedure introduced in \cite{cat1}.
For the convenience of the reader we briefly recall some of 
its key properties.
For a classical map $A \equiv I \mod 2$ and Planck's constant of the
form $h=1/p$, the quantization associates a unitary operator
$\weil_p(A)$ acting on $\HP$ in such a way that the two following
important properties hold: $\weil_p(A)$ only depends on the value of
$A$ modulo $2p$, and the quantization is multiplicative in the sense
that $\weil_p(A B) = \weil_p(A ) \weil_p(B)$ if $A, B$ are two
different classical maps (both congruent to $I$ modulo $2$.)
In fact, $\weil_p(A)$ can be defined via the Weil representation of
$SL_2(\fp)$ acting on $\HP=L^2(\fp)$.  (We abuse notation and also let $A$
denote the image of $A$ in $SL_2(\fp)$ under the reduction modulo $p$
map.)  Since $SL_2(\fp)$ is generated by matrices of the form
\begin{equation}\label{gensodd}
\begin{pmatrix} 1&b\\0 &1 \end{pmatrix},\quad 
\begin{pmatrix} t&0\\0 &t^{-1} \end{pmatrix},\quad  
\begin{pmatrix} 0 &1\\ -1& 0 \end{pmatrix}
\end{equation}
it suffices to specify $\weilp$ on such matrices. 
Let $e(x) = e^{2\pi i x}$, and let $r$ be an integer such that $2r
\equiv 1 \mod p$.  Further, with
$$S_r(a,p) := 
\frac 1{\sqrt{p}} \sum_{x\bmod p} e(\frac {-ra x^2}{p}),
$$
and $\Lambda$  being the nontrivial quadratic character on $\fp^\times$, 
the
action on a state $\phi \in \HP$ is given by
\begin{equation}\label{U(1b01)}
\left( \weilp \left( \begin{pmatrix} 1&b\\ 0&1 \end{pmatrix} \right) \phi
\right) (x) 
=
e(\frac{rbx^2}{p})\phi(x)
\end{equation}
\begin{equation}\label{U(t001/t)}
\left(\weilp \left( \begin{pmatrix} t&0\\ 0&t^{-1} \end{pmatrix} \right)
\phi \right)(x) =
\Lambda(t) \phi(tx)
\end{equation}
\begin{equation}\label{U(w)}
\left( \weilp \left( \begin{pmatrix}0 &1\\ -1&0 \end{pmatrix} \right)
\phi \right) (x) = 
S_r(-1,p) \frac 1{\sqrt{p}}
\sum_{y\bmod p}\phi(y)e(\frac {2rxy}{p}).
\end{equation}
To simplify the notation, we let $\psi(x) = e(rx/p)$; $\psi$ is then a
nontrivial additive character on $\fp$.

As an immediate consequence, we obtain the following description of
how the Weil representation acts on  delta functions:
\begin{lem} 
Let $\delta_i$ be a delta function supported at $i$, i.e.
$\delta_i(x) = 1$ if $x \equiv i \mod p$, and $\delta_i(x) = 0$
otherwise.
Let $M = \begin{pmatrix} a&b\\c&d\end{pmatrix}$ be an element in
$SL_2(\fp)$ such that $c \not \equiv 0 \mod p$. Then
\begin{equation}
\label{e:weil-action}    
\left( \weilp(M)
\delta_i \right) (x) 
= 
\frac{S_r(-1,p)}{\sqrt{p}}
\Lambda(-c) \psi \left(
\frac{a x^2 + di^2 - 2xi}{c} \right)  
\end{equation}
\end{lem}
\begin{proof}
Since $c \not \equiv 0 \mod p$, we can write $M$ as the following
product of generators
\begin{equation}
  \label{eq:weil-general-action}
M = 
\begin{pmatrix} a&b\\c&d\end{pmatrix}
=
\begin{pmatrix}
1 & b_1 \\ 0 & 1
\end{pmatrix}
\begin{pmatrix}
0 & 1 \\ -1 & 0
\end{pmatrix}
\begin{pmatrix}
1 & b_2 \\ 0 & 1
\end{pmatrix}
\begin{pmatrix}
t & 0 \\ 0 & t^{-1}
\end{pmatrix}
\end{equation}
where $t = -c, b_1 = a/c$ and $b_2 = cd$.  Hence
\begin{equation}
\label{eq:3}
(\weilp(M)\delta_i)(x) = 
\frac{S_r(-1,p)}{\sqrt{p}}
 \Lambda(t) \sum_{y} \psi( b_1x^2 + b_2y^2 + 2xy) \delta_i(ty)
\end{equation}
and since the only nonvanishing term in the sum over $y$ is when
$ty=-cy = i$, we find that 
$$
(\weilp(M)\delta_i)(x) = 
\frac{S_r(-1,p)}{\sqrt{p}}
\Lambda(-c) \psi \left(
\frac{a}{c}  x^2 + cd (\frac{i}{c})^2 + 2x\frac{i}{-c} \right) 
$$
$$
=
\frac{S_r(-1,p)}{\sqrt{p}}
\Lambda(-c) \psi \left(
\frac{a x^2 + di^2 - 2xi}{c} \right)  
$$

\end{proof}
\begin{remark}
\label{rem:trace-calc}  
As an immediate consequence of the Lemma we can, up to a phase,
determine the trace of $\weilp(M)$ in many cases:
$$
\left| \trace \left( 
\weilp \left( \begin{pmatrix} a & b \\ c&d\end{pmatrix} \right)
\right) \right|
=1
$$
if $c \neq 0$ and $a+b \neq 2$ since the absolute value of the
Gauss sum $ \sum_{x \mod p} \psi( \alpha x^2)$ equals
$\sqrt{p}$ if $\alpha \neq 0$.  Moreover, from (\ref{U(t001/t)}),
we see that $|\trace(\weilp(-I))| = 1$, and trivially, $\trace(\weilp(I))
= p$.
\end{remark}

\begin{remark}
\label{rem:no-parity}  
We may relax the condition $A \equiv I \mod 2$ by defining $\torus$ in
a different way and using a slightly different quantization procedure:
the image of $A$ in $SL_2(\Z/2p\Z)$ is contained in some cyclic group
of maximal order, say generated by some element $B$.  Define $\torus$
to be this cyclic subgroup.  The Weil representation then gives a
quantization $\weil_p(B)$ which is well defined up to a choice of
scalar, and this scalar can be chosen so that multiplicativity
holds on $\torus$.

\end{remark}

\subsection{Hecke operators and eigenfunctions}
\label{sec:hecke-eigenfns}
Let $p$ be a fixed inert prime (i.e., the characteristic polynomial of
$A$ remains irreducible modulo $p$.)  With 
$$
\torus := \{ B \in SL_2(\fp) : AB  \equiv BA \mod p\},
$$
the Hecke operators are then given by $\{ \weil_p(B) : B \in
\torus\}$.  
(In \cite{cat1}, the Hecke operators were defined in a
somewhat different way: a subgroup of the centralizer of $A$ modulo
$2N$ was identified with the norm one elements of $\OO/(2N\OO)$ (with
certain parity conditions) where $\OO$ is an order contained in
$\Q(\epsilon)$ and $\epsilon$ is an eigenvalue of $A$.  However, in the
case of $N=p$ an inert prime it is straightforward to verify that
these notions are the same.)

A Hecke eigenfunction is then a state $\phi_\nu \in \HP$ such that
$$
\weil_p(B) \phi_\nu = \nu(B) \phi_\nu  \quad \forall B \in \torus
$$
where $\nu : \torus \to \C^\times$ is a character of $\torus$.  (Note
that the eigenspaces of the Hecke operators are paremetrized by
characters of $\torus$.)  Unless otherwise noted, all Hecke
eigenfunctions will be normalized so that $\norm{\phi_\nu} = 1$.

Our goal is to express Hecke eigenfunctions in terms of exponential
sums in one variable, and in order to do this we need a geometric
parametrization of $\torus$.  Since we assume that $p$ is inert, there
exist some $M=
\begin{pmatrix}X & Y\\Z & W\end{pmatrix} \in SL_2(\fp)$ so that
\begin{equation}
\label{e:conjugate-torus}
\torus = M
\left\{
\begin{pmatrix}a &
  bD\\b&a
\end{pmatrix} : a^2 - Db^2 \equiv 1 \mod p \right\}
M^{-1}
\end{equation}
where $D$ is not a square in $\fp$.  
%
We note that the solutions to $a^2 - Db^2 =1$ can be
parametrized by letting 
$$
(a,b) = \left( \frac{1+Dt^2}{1-Dt^2}, \frac{2t}{1-Dt^2} \right)
$$
where $t$ takes values in $P^1(\fp)$.
We next determine how Hecke operators act on
delta functions.
\begin{lem}
Let $M=\begin{pmatrix}X & Y\\Z & W\end{pmatrix}$ and $B=\begin{pmatrix}a
  &  bD\\b&a\end{pmatrix}$ be elements in $ SL_2(\fp)$. If $b \neq 0$,
  then
\begin{multline}
\label{e:simpleaction}
\left( \weilp(MBM^{-1})
\delta_i \right) (x) 
=
\Lambda(W^2-DZ^2)
\psi \left(
\frac{(YW-DXZ)(x^2-i^2)  }
{(W^2-DZ^2)} 
\right) \times
\\
\times
S_r(-1,p)
\Lambda(-b) 
\psi \left(
\frac{
a (x^2 +  i^2) - 2xi}{b(W^2-DZ^2)} \right)
\end{multline}
\end{lem}

\begin{proof}
Since
\begin{multline*}
MBM^{-1} =
\begin{pmatrix}X & Y\\Z &W\end{pmatrix}
\begin{pmatrix}a &bD\\b&a\end{pmatrix}
\begin{pmatrix}W &-Y\\ -Z&X\end{pmatrix}
\\=
\begin{pmatrix}a(XW-YZ)+b(YW-DXZ) & *       \\
b(W^2-DZ^2) & a(XW-YZ)+b(DXZ-YW)\end{pmatrix}
\\
=
\begin{pmatrix}a+b(YW-DXZ) & *       \\
b(W^2-DZ^2) & a+b(DXZ-YW) \end{pmatrix},
\end{multline*}
equation (\ref{e:weil-action}) gives  that\footnote{Note that
  $b(W^2-DZ^2) \neq 0$ since $b \neq 0$ and 
$W^2-DZ^2 \neq 0$ since $D$ is not a square in $\fp$.}
\begin{multline*}
\left( \weilp(MBM^{-1})
\delta_i \right) (x) = 
S_r(-1,p)
\Lambda(-b(W^2-DZ^2))  \times
\\ 
\times
\psi \left(
\frac{
(a+b(YW-DXZ)) x^2 + (a+b(DXZ-YW) )i^2 - 2xi}{b(W^2-DZ^2)} \right)
=
\\
=
S_r(-1,p)
\Lambda(W^2-DZ^2)
\Lambda(-b) 
\psi \left(
\frac{(YW-DXZ)(x^2-i^2)  }
{(W^2-DZ^2)} 
\right)
\times \\ \times
\psi \left(
\frac{
a (x^2 + i^2) - 2xi}{b(W^2-DZ^2)} \right)
\end{multline*}
\end{proof}
%

\section{Proof of Theorem~\ref{thm:main}}

\subsection{Hecke eigenfunctions via projections}
\label{subsec:constr-eigenf}
Given a character $\nu$ on $\torus$, let 
$$V_\nu := \{ \phi \in L^2(\fp) :
\weilp(B) \phi = \nu(B) \phi\ \ \forall B \in \torus\}.
$$ 
Since $p$ is odd and inert in $K$, the Weil representation
restricted to $\torus$ 
is multiplicity free and hence
$\operatorname{dim}(V_\nu) \leq 1$.
For a proof of this, see Proposition~3 in \cite{moen-dual-pair-u1-u1}.
Alternatively, 
we might argue as follows: let $d_\nu$ be the dimension of $V_\nu$.
Then $\sum_\nu d_\nu = \dim(L^2(\fp)) = p$, and by the character
formula for group representations, $\sum_\nu d_\nu^2 =
|\torus|^{-1} \sum_{B \in \torus} |\trace(\weilp(B))|^2$, which
by Remark~\ref{rem:trace-calc} equals $|\torus|^{-1}(1 \cdot p^2 +
(p-1)\cdot 1 + 1 \cdot 1) = p(p+1)/|\torus|$.  Since $|\torus| = p+1$
in the inert case, we find that $\sum_\nu d_\nu = p = \sum_\nu d_\nu^2 $
which implies that $d_\nu \leq 1$.
%
In fact, since $\dim(L^2(\fp))=p$ has dimension $p$, we note that
$\operatorname{dim}(V_\nu) =1$ for all but one character $\nu$
of $\torus$.  In what follows we let
$\nu$ be a fixed character of 
$\torus$ such that $\operatorname{dim}(V_\nu) = 1$.

To construct the Hecke eigenfunction corresponding to $\nu$ we define 
a projection operator $\PP : L^2(\fp) \to V_{\nu}$ by letting
$$
\PP f := 
\frac{1}{|\torus|}
\sum_{B \in \torus} \overline{\nu}(B) \weilp(B)f
$$
Clearly,  $\PP f$ is a Hecke eigenfunction;  the
main difficulty is 
to control the $L^2$-norm of $\PP f$.
%
However, if $f$ is a delta function, the $L^2$-norm of the projection
can be expressed in a simple manner.
\begin{lem} 
\label{lem:l2-norm}
We have
$$
\norm{\PP \delta_i}^2 = 
\frac{(\PP \delta_i)(i)}{p}
$$
\end{lem}
\begin{proof} Since $\PP$ is an orthogonal projection\footnote{
The inner product defined by (\ref{eq:inner-product}) is
    invariant by the action of $\torus$, hence the Hilbert space of
    states decomposes into an orthogonal  sum of
    $\torus$-invariant subspaces.}, $\PP^2=\PP$ and $\PP$ is
  self adjoint.  Thus
$$
\norm{\PP \delta_i}^2 = 
\langle \PP \delta_i, \PP \delta_i \rangle =
\langle \PP^2 \delta_i, \delta_i \rangle =
$$
$$ =
\langle \PP \delta_i, \delta_i \rangle 
=
\frac{1}{p} \sum_{x \in \fp} 
(\PP \delta_i)(x) \overline{\delta_i(x)}
=
\frac{(\PP \delta_i)(i)}{p}
$$
\end{proof}

\subsection{Hecke eigenfunctions and exponential sums}
\label{sec:exponential-sums}

We can now express the Hecke eigenfunctions in terms of exponential
sums.  
Since we may regard $\torus$ as $\fpn$, the group of norm one
elements in $\fpt^\times$, and $\fpt = \fp[\sqrt{D}]$, 
Hilbert's theorem 90 gives us a parametrization $P^1(\fp) \to \fpn$
via the map $t \to \frac{1+t\sqrt{D}}{1-t\sqrt{D}}$.

\begin{prop}  
\label{prop:projection-as-exp-sum}
Given a character $\nu : F_{p^2}^1 \to \C^\times$ on
  the group of norm one elements in $F_{p^2}$, define an exponential
  sum 
$$
\Exp(i,x) := 
\sum_{t \neq 0}
\overline{\nu} \left(  
\frac{1+t\sqrt{D}}{1-t\sqrt{D}}
\right)
\Lambda \left( \frac{-2t}{1-Dt^2} \right)
\psi_{\CA} \left( \frac{(x-i)^2}{2t} + Dt(x+i)^2   \right)
$$
where $\CA = (W^2-DZ^2)^{-1}$ and
$\psi_\CA : \fp \to \C^\times$ is a nontrivial additive character
defined by $\psi_\CA(x) 
= \psi(\CA x)$.   Putting 
$$\alpha(i,x) = \Lambda(W^2-DZ^2)
\psi \left( \frac{(YW-DXZ)(x^2-i^2)}{W^2-DZ^2} \right),$$
we then have 
\begin{multline}
\label{eq:1}
(\PP \delta_i)(x)
= \\ =
\frac{1}{p+1} \left(
\delta_i(x) +
\Lambda(-1) \delta_i(-x)
+
\alpha(i,x) S_r(-1,p) 
\frac{\Exp(i,x)}{\sqrt{p}}
\right)  
\end{multline}
\end{prop}
\begin{proof}
We first note that $\pm I \in \torus$ corresponds to $t = 0$ or $t
= \infty$ in the parametrization $P^1(\fp) \to 
\torus$.  Since  $\weilp(I) \delta_i =\delta_i$ and $(\weilp(-I) \delta_i)(x) =
\Lambda(-1) \delta_i(-x)$  we find that 
$$
(\PP \delta_i)(x) = 
\frac{1}{p+1} \left(
\delta_i(x) +
\Lambda(-1) \delta_i(-x)
+
\sum_{B \in \torus, B \neq \pm I} \overline{\nu}(B) 
\left( \weilp(B) \delta_i \right)(x)
\right) 
$$
If $B \in \torus$ and  $B \neq \pm I$, then
$$
B = M
\begin{pmatrix}a &
  bD\\b&a
\end{pmatrix} 
M^{-1}
$$
where $ a^2 - Db^2 \equiv 1 \mod p$ and $b \not \equiv 0 \mod p$.  
These solutions are all of the form
$(a,b) = \left( \frac{1+Dt^2}{1-Dt^2}, \frac{2t}{1-Dt^2} \right)$
for $t \in \fp^\times$.  
From (\ref{e:simpleaction}) we obtain
\begin{multline*}
\overline{\nu}(B)
\left( \weilp(B) \delta_i \right)(x)
\\ = 
\overline{\nu} \left(  
\frac{1+t\sqrt{D}}{1-t\sqrt{D}}
\right)
\alpha(i,x) 
\Lambda(-b) 
S_r(-1,p)
\psi \left(
\frac{
a (x^2 + i^2) - 2xi}{b(W^2-DZ^2)} \right)
\end{multline*}
and since 
$$
\frac{a(x^2+i^2)-2xi}{b} = 
\frac{(1+Dt^2)(x^2+i^2)-2xi(1-Dt^2)}{2t} 
$$
$$= 
\frac{Dt(x+i)^2}{2} + \frac{(x-i)^2}{2t}.
$$
the result follows.
\end{proof}

{\em Remark:}  If $x=i$ then $\alpha(i,x)= \Lambda(W^2-DZ^2) = \pm
1$, and it can be shown that $S_r(-1,p) \Exp(i,i)$ is in fact real
valued.

Using the Riemann hypothesis for curves, we will now bound the
exponential sums $\Exp(i,x)$.
\begin{prop}
\label{prop:exp-sum-bound}
If $x \neq \pm i$, then 
$$
|\Exp(i,x)| \leq 4 \sqrt{p},
$$
and if $x= \pm i$, then
$$
|\Exp(i,x)| \leq 3 \sqrt{p}.
$$
Moreover,
$$
|\Exp(0,0)| \leq 2 \sqrt{p}.
$$
\end{prop}

\begin{proof}
Let $w$ be the place corresponding to the
field extension 
$\fpt/\fp$.  
%
By chapter 6 of \cite{winnie-li-book} (see Theorems 4, 6, and exercise
3) there exists idele class characters $\tilde{\nu},
\tilde{\psi},\tilde{\Lambda}$ such that
$$
\tilde{\nu}(\pi_v) =
\overline{\nu} \left(  
\frac{1+t\sqrt{D}}{1-t\sqrt{D}}
\right) 
$$
$$
\tilde{\psi}(\pi_v) =
\psi \left( \frac{(x-i)^2}{2t} + Dt(x+i)^2   \right) 
$$
$$
\tilde{\Lambda}(\pi_v) = \Lambda \left( \frac{-2t}{1-Dt^2} \right)
$$
for all degree one places in $\fp[X]$ of the form
$\pi_v = (X+t)$, $t \in \fp^\times$.
Moreover, the conductors are as follows: 
$$\cond(\tilde{\nu}) = (w),$$
$$\cond(\tilde{\Lambda})=(0)+(w)+(\infty),$$
and 
$$
\cond(\tilde{\psi}) =
\begin{cases}
2(0)+2(\infty) &\text{ if } x \neq \pm i \\
2(\infty) &\text{ if } x=i, \\
2(0) &\text{ if } x=-i, \\
0 &\text{ if } x=i=0. \text{ (I.e., $\tilde{\psi}$ is trivial.)} \\
\end{cases}
$$
Letting $\chi = \tilde{\nu} \cdot \tilde{\psi} \cdot \tilde{\Lambda}$
we have 
$$
\Exp(i,x) = 
\sum_{\substack{v : \deg(v)=1  \\v \text{ unramified}}}
\chi(\pi_v)
$$
and Corollary 3 in 6.1 of \cite{winnie-li-book} then gives 
$$
\left|
\sum_{\substack{v : \deg(v)=1  \\v \text{ unramified}}}
\chi(\pi_v)
\right| 
\leq 
\sqrt{p} (m-2)
$$
where $m$ is the degree of the conductor of $\chi$.
Now, if $i \neq \pm x$ then  $\cond(\chi) = (w)+2(0)+2(\infty)$,
which has degree $2+2+2 = 6$, and hence  
$$
|\Exp(i,x)| \leq 4 \sqrt{p} 
$$
If $x=i$, then $\cond(\chi) = (w)+(0)+2(\infty)$, and if $x=-i$, then 
$\cond(\chi) = (w)+2(0)+(\infty)$. In either case, the degree is $5$, 
hence 
$$
|\Exp(i,\pm i)| \leq 3 \sqrt{p}.
$$
Finally, a similar argument gives that $|\Exp(0,0)| \leq  2\sqrt{p}$.

\end{proof}

\subsection{Conclusion}
Let $f_{\nu,i} = \PP \delta_i$.  Since $f_{\nu,i}$ is always a scalar
multiple (possibly zero) 
of an $L^2$-normalized eigenfunction
$\phi_{\nu}$, we may find functions $g_\nu, h_\nu$ such that
$$
f_{\nu,i}(x) =
g_\nu(i) h_\nu(x)
$$
where $h_\nu$ is a scalar multiple of $\phi_v$ normalized so that
$
\sum_x |h_\nu(x)|^2  = 1.
$
Since
$$
f_{\nu,i}(x) = 
p \langle P_\nu \delta_i, \delta_x \rangle =  
p \langle  \delta_i, P_\nu \delta_x \rangle =  
p \overline{\langle  P_\nu \delta_x,  \delta_i \rangle} =  
\overline{f_{\nu,x}(i)} 
$$
we find that
$$
g_\nu(i) h_\nu(x) = \overline{g_\nu(x)}  \overline{h_\nu(i)}
$$
for all $i,x$.  Now, if $\PP \neq 0$ then $\PP \delta_i \neq
0$ for some $i$, hence $g(i) \neq 0$ for some $i$ and we find that
$$
h_\nu(x) =    \overline{g_\nu(x)} \cdot \frac{\overline{h_\nu(i)}} {g_\nu(i)}
$$
for all $x$.
In order to determine $\overline{h_\nu(i)}/g_\nu(i)$ we argue
as follows: since we have chosen $h_\nu$ so that $\sum_x |h_\nu(x)|^2  = 1$,
$$
\norm{f_{\nu,i}}^2 = 
\frac{1}{p}
\sum_x |f_{\nu,i}(x)|^2 
= 
\frac{|g_\nu(i)|^2}{p}
\sum_x |h_\nu(x)|^2 
=
\frac{|g_\nu(i)|^2}{p},
$$
which, using Lemma~\ref{lem:l2-norm}, also can be written as
$$
\norm{f_{\nu,i}}^2 
=
\norm{\PP \delta_i}^2
=
\frac{(\PP \delta_i)(i)}{p}
=
\frac{f_{\nu,i}(i)}{p}
=
\frac{g_\nu(i) h_\nu(i)}{p}.
$$
Thus $h_\nu(i) = \overline{g_\nu(i)}$
and we find that 
$h_\nu(x) = \overline{g_\nu(x)}$ for {\em all} $x$.
Hence
$$
|h_\nu(x)|^2 = |h_\nu(x) g_\nu(x)| = |f_{\nu,x}(x)|
=
|(\PP \delta_x)(x)|,
$$
and by Proposition~\ref{prop:projection-as-exp-sum} we find that
$$
|h_\nu(x)|^2
=
\begin{cases}
\left| \frac{1}{p+1} \left(
1
+
\alpha(x,x) S_r(-1,p) 
\frac{\Exp(x,x)}{\sqrt{p}}
\right)
\right| & \text{ if $x \neq 0$,}
\\
\left| \frac{1}{p+1} \left(
1 + \Lambda(-1) +
\alpha(0,0) S_r(-1,p) 
\frac{\Exp(0,0)}{\sqrt{p}}
\right)
\right|
& \text{ if $x = 0$.}
\end{cases}
$$
By Proposition~\ref{prop:exp-sum-bound}, $|\Exp(x,x)| \leq 3 \sqrt{p}$
if $x \neq 0$, and $|\Exp(0,0)| \leq 2 \sqrt{p}$, 
hence
$$
|\phi_{\nu}(x)| 
=
p^{1/2}|h_\nu(x)| 
\leq
2 \sqrt{\frac{p}{p+1}}
$$
since $|S_r(-1,p)|=1$, and $|\alpha(x,x)|=1$ for all $x$.

\section{Proof of Theorem~\ref{thm:summary}}
There are three different types of primes that needs to be considered:
inert, split and ramified primes.  Theorem~\ref{thm:main} gives the
inert case (note that $A$ cannot be upper triangular modulo $p$ if $p$
is inert), and the case of $A$ split and not upper triangular is
Theorem~2 of \cite{catsup}.  What remains to be done is to treat the
ramified case, and the case of $p$ split and $A$ upper triangular
modulo $p$.

\subsection{The ramified case}
If the characteristic polynomial of $A$ has a double root modulo $p$
then $p$ is said to be ramified.  In this case, 
$$
A = M \begin{pmatrix}\pm 1 & s\\ 0 & \pm 1\end{pmatrix} M^{-1}
$$
for some $M \in SL_2(\fp)$ and $s \in \fp$.  Moreover, the norm one
elements used to define the Hecke operators correspond to matrices
conjugate to upper triangular matrices with $\pm 1$ on the diagonal,
i.e., 
$$
\torus = M \left\{ 
\begin{pmatrix}\pm 1 & t\\ 0 & \pm 1\end{pmatrix}
: t \in \fp \right\} M^{-1}.
$$
The normalized Hecke eigenfunctions are then given by
$$
\phi_i^\pm = \sqrt{p/2} \cdot \weil_p(M)( \delta_i \pm \delta_{-i})
$$
for $0 \leq i \leq (p-1)/2$ and 
$$
\phi_0 = \sqrt{p} \cdot \weil_p(M) \delta_0
$$

We note that $A$ is upper triangular if and only if $M$ is upper
triangular.  In the upper triangular case, $\weil_p(M)$ acts via
multiplication by scalars (of absolute value $1$), and by permuting
the arguments (see equations (\ref{U(1b01)}) and (\ref{U(t001/t)})), 
hence 
$$
\supnorm{\phi_i^\pm} = \sqrt{p/2} \supnorm{ \delta_i \pm \delta_{-i}}
= \sqrt{p/2}
$$
and 
$$
\supnorm{\phi_0} = \sqrt{p} \supnorm{ \delta_0} = \sqrt{p}
$$

If $A$ is not upper triangular, then $M = \begin{pmatrix} a & b \\
  c & d\end{pmatrix}$ for some $c \not \equiv 0 \mod p$, and 
equation~\ref{e:weil-action} then gives that
$$
|(\weil_p(M)(\delta_i \pm \delta_{-i})(x)| =
\frac{1}{\sqrt{p}} 
| \psi(\frac{ax^2+di^2 -2xi}{c}) + \psi(\frac{ax^2+di^2 +2xi}{c})|
$$
and thus 
$$
\supnorm{\phi_i^\pm} = \sqrt{2},
$$
(for instance, take $x = 0$) and a similar calculation gives that 
$$
\supnorm{\phi_0} = 1.
$$

\subsection{The upper triangular split case}
We assume that $A$ is upper triangular and conjugate to a diagonal
matrix, i.e., 
$$
A = M  \begin{pmatrix}a & 0\\ 0 & a^{-1}\end{pmatrix}M^{-1}
$$
where $M$ is upper triangular.  Arguing as in section 4.1 of
\cite{catsup}, we find that a basis of normalized Hecke eigenfunctions
are given by
$$
\phi_\chi = \sqrt{\frac{p}{p-1}} \weil_p(M) \chi
$$
where $\chi : \fp \to \C^\times$ ranges over all multiplicative
characters of $\fp^\times$ (extended to $\fp$ by letting $\chi(0)=0$)
and
$$
\phi_0 = \sqrt{p} \cdot  \weil_p(M) \delta_0.
$$
Since  $M$ is upper triangular, 
$\weil_p(M)$ acts via
multiplication by scalars of absolute value $1$ and by permuting
the arguments, hence
$$
\supnorm{\phi_\chi} = \sqrt{\frac{p}{p-1}}
$$
and 
$$
\supnorm{\phi_0} = \sqrt{p}
$$
{\em Remark:} Since one of the eigenspaces have dimension two, we can
find a character $\chi$ and choose scalars $\alpha,\beta$ so that
$\phi = \alpha \phi_0 + \beta 
\phi_\chi$ is a normalized Hecke eigenfunction with supremum norm
equal to $\sqrt{p^2/(p-1)}$.

\section{Proof of Theorem~\ref{thm:square-free-bound}}
If $N$ is square free, say with prime factorization $N = \prod_{i=1}^k
p_i$, then $\UN(A)$ can be expressed as a tensor product (cf.
\cite{cat1}, section~4.1) of the form $ \UN(A) = U_{p_1}(A)
\otimes U_{p_2}(A) \otimes \cdots \otimes U_{p_k}(A)$.
Thus any Hecke eigenfunction $\phi \in \HN$ can be written as a
product of Hecke functions $\phi_i \in H_{p_i}$, i.e., $\phi(x)
= \prod_{i=1}^k \phi_i(x_i)$ where each $x_i$ is the image of $x$
under the projection from $\Z/N\Z$ to $\Z/p_i\Z$.  Since $A$ cannot
be upper triangular modulo $p$ for more than a finite number of
primes, Theorem~\ref{thm:summary}  gives that 
$
\supnorm{\phi} = \prod_{i=1}^k \supnorm{\phi_i} \ll 2^k \ll_\epsilon N^\epsilon.
$


\begin{thebibliography}{10}


\bibitem{aurich-backer-schubert-taglieber-maximum-norms}
R.~Aurich, A.~B{\"a}cker, R.~Schubert, and M.~Taglieber.
\newblock Maximum norms of chaotic quantum eigenstates and random waves.
\newblock {\em Phys. D}, 129(1-2):1--14, 1999.

\bibitem{aurich-steiner-statistical-properties-highly-excited}
R.~Aurich and F.~Steiner.
\newblock Statistical properties of highly excited quantum eigenstates of a
  strongly chaotic system.
\newblock {\em Phys. D}, 64(1-3):185--214, 1993.

\bibitem{berard}
P.~H. B{\'e}rard.
\newblock On the wave equation on a compact {R}iemannian manifold without
  conjugate points.
\newblock {\em Math. Z.}, 155(3):249--276, 1977.

\bibitem{berry-regular-irregular-semiclassical}
M.~V. Berry.
\newblock Regular and irregular semiclassical wavefunctions.
\newblock {\em J. Phys. A}, 10(12):2083--2091, 1977.

\bibitem{DE}
M.~Degli~Esposti.
\newblock Quantization of the orientation preserving automorphisms of the
  torus.
\newblock {\em Ann. Inst. H. Poincar\'e Phys. Th\'eor.}, 58(3):323--341, 1993.

\bibitem{DEGI}
M.~Degli~Esposti, S.~Graffi, and S.~Isola.
\newblock Classical limit of the quantized hyperbolic toral automorphisms.
\newblock {\em Comm. Math. Phys.}, 167(3):471--507, 1995.

\bibitem{gurevich-hadani-hecke-sup-norms}
S.~Gurevich and R.~Hadani.
\newblock Heisenberg {R}ealizations, {E}igenfunction and {P}roof of the
  {K}urlberg-{R}udnick {S}upremum {C}onjecture.
\newblock {\em Preprint (preliminary version)}.

\bibitem{gurevich-hadani-hanny-berry-equivariant}
S.~Gurevich and R.~Hadani.
\newblock On {H}annay-{B}erry {E}quivariant {Q}uantization of the {T}orus.
\newblock {\em Preprint}, 2002.

\bibitem{BH}
J.~H. Hannay and M.~V. Berry.
\newblock Quantization of linear maps on a torus-{F}resnel diffraction by a
  periodic grating.
\newblock {\em Phys. D}, 1(3):267--290, 1980.

\bibitem{hejhal-rackner-gaussian-value-distribution}
D.~A. Hejhal and B.~N. Rackner.
\newblock On the topography of {M}aass waveforms for ${P}{S}{L}(2,\bold {Z})$.
\newblock {\em Experiment. Math.}, 1(4):275--305, 1992.

\bibitem{iwaniec-sarnak-supnorms}
H.~Iwaniec and P.~Sarnak.
\newblock {$L\sp \infty$} norms of eigenfunctions of arithmetic surfaces.
\newblock {\em Ann. of Math. (2)}, 141(2):301--320, 1995.

\bibitem{katz-kr-sums}
N.~M. Katz.
\newblock Sato-{T}ate equidistribution of {K}urlberg-{R}udnick sums.
\newblock {\em Internat. Math. Res. Notices}, (14):711--728, 2001.

\bibitem{klimek-leniewski-maitra-rubin}
S.~Klimek, A.~Le{\'s}niewski, N.~Maitra, and R.~Rubin.
\newblock Ergodic properties of quantized toral automorphisms.
\newblock {\em J. Math. Phys.}, 38(1):67--83, 1997.

\bibitem{Knabe}
S.~Knabe.
\newblock On the quantisation of {A}rnold's cat.
\newblock {\em J. Phys. A}, 23(11):2013--2025, 1990.

\bibitem{cat1}
P.~Kurlberg and Z.~Rudnick.
\newblock Hecke theory and equidistribution for the quantization of linear maps
  of the torus.
\newblock {\em Duke Math. J.}, 103(1):47--77, 2000.

\bibitem{catsup}
P.~Kurlberg and Z.~Rudnick.
\newblock Value distribution for eigenfunctions of desymmetrized quantum maps.
\newblock {\em Int. Math. Res. Not.}, 2001(18):985--1002, 2001.

\bibitem{winnie-li-book}
W.~C.~W. Li.
\newblock {\em Number theory with applications}.
\newblock World Scientific Publishing Co. Inc., River Edge, NJ, 1996.

\bibitem{mezzadri-multiplicative-cat-maps}
F.~Mezzadri.
\newblock On the multiplicativity of quantum cat maps.
\newblock {\em Nonlinearity}, 15(3):905--922, 2002.

\bibitem{moen-dual-pair-u1-u1}
C.~Moen.
\newblock The dual pair {$({\rm U}(1),{\rm U}(1))$} over a {$p$}-adic field.
\newblock {\em Pacific J. Math.}, 158(2):365--386, 1993.

\bibitem{jakobson-survey}
N.~Nadirashvili, D.~Tot, and D.~Yakobson.
\newblock Geometric properties of eigenfunctions.
\newblock {\em Uspekhi Mat. Nauk}, 56(6(342)):67--88, 2001.

\bibitem{olofsson-big-catmap-supnorm}
R.~Olofsson.
\newblock Large supremum norms for hecke eigenfunctions of quantized cat maps.
\newblock {\em Preprint}.

\bibitem{rudnick-sarnak-eigenstates-arithmetic-hyperbolic-mflds}
Z.~Rudnick and P.~Sarnak.
\newblock The behaviour of eigenstates of arithmetic hyperbolic manifolds.
\newblock {\em Comm. Math. Phys.}, 161(1):195--213, 1994.

\bibitem{sarnak-arithmetic-qc}
P.~Sarnak.
\newblock Arithmetic quantum chaos.
\newblock In {\em The Schur lectures (1992) (Tel Aviv)}, pages 183--236.
  Bar-Ilan Univ., Ramat Gan, 1995.

\bibitem{sarnak-spectra-of-hyperbolic-surfaces}
P.~Sarnak.
\newblock Spectra of hyperbolic surfaces.
\newblock {\em Bull. Amer. Math. Soc. (N.S.)}, 40(4):441--478 (electronic),
  2003.

\bibitem{sogge-sup-norm}
C.~D. Sogge.
\newblock Concerning the {$L\sp p$} norm of spectral clusters for second-order
  elliptic operators on compact manifolds.
\newblock {\em J. Funct. Anal.}, 77(1):123--138, 1988.

\bibitem{sogge-zelditch-max-eigenfunction-growth}
C.~D. Sogge and S.~Zelditch.
\newblock Riemannian manifolds with maximal eigenfunction growth.
\newblock {\em Duke Math. J.}, 114(3):387--437, 2002.

\bibitem{toth-zelditch-bounded-eigenfunctions}
J.~A. Toth and S.~Zelditch.
\newblock Riemannian manifolds with uniformly bounded eigenfunctions.
\newblock {\em Duke Math. J.}, 111(1):97--132, 2002.

\bibitem{toth-zelditch-revisited}
J.~A. Toth and S.~Zelditch.
\newblock Norms of modes and quasi-modes revisited.
\newblock In {\em Harmonic analysis at Mount Holyoke (South Hadley, MA, 2001)},
  volume 320 of {\em Contemp. Math.}, pages 435--458. Amer. Math. Soc.,
  Providence, RI, 2003.

\bibitem{zelditch-quantum-ergodicity-and-mixing}
S.~Zelditch.
\newblock Quantum ergodicity and mixing, to appear in {E}ncyclopedia of
  {M}athematical {P}hysics.

\bibitem{Zelditch97}
S.~Zelditch.
\newblock Index and dynamics of quantized contact transformations.
\newblock {\em Ann. Inst. Fourier (Grenoble)}, 47(1):305--363, 1997.


\end{thebibliography}

\bibliographystyle{abbrv} 


\end{document}